\documentclass[a4paper, 11pt]{article}

\usepackage{rotating}
\usepackage{amsmath}
\usepackage{amssymb}
\usepackage{latexsym}
\usepackage{enumerate}
\usepackage{graphics}
\usepackage{amsmath}
\usepackage{amssymb}
\usepackage{subfig}
\usepackage{graphicx}
\usepackage{todonotes}
\usepackage{amsmath}
\usepackage{amssymb}
\usepackage{graphicx}

\usepackage[english]{babel}
\usepackage{lineno}
\usepackage{color}

\usepackage[utf8]{inputenc}
\usepackage[english]{babel}
\allowdisplaybreaks

\newtheorem{theo}{Theorem}

\newtheorem{lem}{Lemma}

\newtheorem{prop}{Proposition}

\newtheorem{defi}{Definition}

\newtheorem{cor}{Corollary}

\newcommand{\la}{\langle}
\newcommand{\ra}{\rangle}

\begin{document}
%\linenumbers*[1]

\title{On bounding exact models of epidemic spread on networks}

\author{P\'eter L. Simon$^{1,2,\ast}$, Istv\'an Z. Kiss$^{3}$}

\maketitle

\begin{center}
$^1$ Institute of Mathematics, E\"otv\"os Loránd University Budapest, Hungary \\
$^2$ Numerical Analysis and Large Networks Research Group,\\ Hungarian Academy of Sciences, Hungary\\
$^3$ School of Mathematical and Physical Sciences, Department of Mathematics, University of Sussex, Falmer, Brighton BN1 9QH, UK
\end{center}

\vspace{1cm}

\begin{abstract}
%Modelling the spread of epidemics on networks has led to a myriad of mean-filed models such as heterogenous degree-based, pairwise, effective-degree, edge-based compartmental and N-Intertwined Mean-Field Approximation (NIMFA). All these deterministic models aim to approximate the underlying exact stochastic process by encoding the structure of the network and the properties of the transmission process. Except a few cases, mean-filed models are validated heuristically by comparing outputs from these to results based on explicit stochastic network simulations, and comparisons are usually only feasible for a limited number of different network types and parameter combinations.
In this paper we use comparison theorems from classical ODE theory in order to rigorously show that the N-Intertwined Mean-Field Approximation (NIMFA) model provides an upper estimate on the exact stochastic process. The proof of the results relies on the observation that the epidemic process is negatively correlated (in the sense that the probability of an edge being in the susceptible-infected state is smaller than the product of the probabilities of the nodes being in the susceptible and infected states, respectively), which we also prove rigorously. The results in the paper hold for arbitrary weighted and directed networks. We cast the results in a more general framework where
alternative closures, other than that assuming the independence of nodes connected by an edge, are possible and provide a succinct summary of the stability analysis of the resulting more general mean-field models.
%Here, wea rigorous results is presented some rigorous results are presented based on two different techniques. The first uses comparison theorems from classical ODE theory,
%while the second starts from the forward Kolmogorov equations and manipulates the moments favourably to bound the absolute difference between the prevalence of the exact and mean-field models. The first approach is more general and it applies for arbitrary networks with the second being limited to fully connected networks.
\end{abstract}


\vspace{1cm}


\vspace{1cm}
\begin{flushleft}
$\ast$ corresponding author\\
email: simonp@cs.elte.hu\\
\end{flushleft}

\newpage

%%%%%%%%%%%%%%%%%%%%%%%%%%%%%%%%%%%%%%%%%%%%%%%%%%%%%%%%%%%%%%%%%%%%%%%%%%%%%%%%%%%%%%%%%%
\section{Introduction}
%%%%%%%%%%%%%%%%%%%%%%%%%%%%%%%%%%%%%%%%%%%%%%%%%%%%%%%%%%%%%%%%%%%%%%%%%%%%%%%%%%%%%%%%%%
Modelling transmission processes on networks, such as epidemics and rumours, has led to many mathematical challenges~\cite{pastor2014epidemic,kiss2016mathematics}.
This is mainly due to the high dimensionality of the exact model, which is often a continuous time Markov chain
where the size of the state space scales as $m^N$, where $m$ is the number of states a node can be in and $N$ is the number of
nodes in the network \cite{simon2011JMB,taylor2012JMB}. While theoretically the master equations can be given, their rigorous analysis is out of reach due to the high dimensionality.
One approach to deal with this challenge is to consider some `clever' averaging, at node or at population level, and proceed to derive evolution equations
for some newly defined average quantities. These however, more often than not, depend on other new average quantities which are of higher order, e.g.
the expected number of nodes in a certain state usually depends on the expected number of links/edges in certain states. As a rule of thumb,
the dependency between moments is broken by making some closure assumptions where higher order moments are approximated by lower-order ones. This then leads to
a low-dimensional system of ordinary differential equations (integro-differential and delay differential equations are also possible) or mean-field model.

The approach above has led to a myriad of mean-field models for $SIS$ (susceptible-infected-susceptible) and $SIR$ (susceptible-infected-recovered) epidemics which are able to capture the average behaviour of epidemics on certain network types (e.g. tree networks and networks build using the configuration model). The most well-known such models are: (a) heterogenous degree-based~\cite{SatorrasVespignani,boguna2002epidemic}, (b) pairwise~\cite{rand1999,Keeling1999,vanBaalen2000,simon2011JMB,taylor2012JMB,kiss2015generalization}, (c) effective-degree~\cite{lindquist2011effective}, (d) edge-based compartmental~\cite{miller2012edge,miller2013model}, (e) pair-based~\cite{Sharkey2011,sharkey2013exact,kiss2015exact} and (f) N-Intertwined Mean-Field Approximation (NIMFA)~\cite{van2009virus,PietComputingSpringer}. When such models perform well, i.e. output from these agrees well with results from the exact or simulation model, one can proceed to analyse them and to derive analytical results concerning the epidemic threshold and final size, or quasi-equilibrium for $SIS$ epidemics. Such explicit relations between network characteristics and spreading dynamics allow us to better understand how these factors interact and will ultimately lead to more efficient prevention and control measures.

In many cases, mean-field models are validated by simply comparing the results of explicit stochastic network simulation (which stands for the exact model) to output from mean-field models.
Such tests are usually performed for a limited number of network types and combination of parameter values. While such heuristic arguments are useful, it is desirable that where possible the difference between the exact and mean-field models is rigorously established using sound mathematical arguments. This endeavour has already led to results proving that under some mild conditions on the degree distribution the edge-based compartmental model is exact for $SIR$ epidemics in the limit of large networks built based on the configuration model~\cite{MolloyReed}. In the case of $SIS$ propagation on a complete graph or on a regular random network the model is a density dependent Markov chain and functional analytic tools can be used to prove that the difference between the output of the mean-field and the exact system scales as $1/N$ for large system size $N$ \cite{batkai2011differential,kunszenti2016mean}. Besides estimates on the difference, upper and lower bounds have also been derived for the prevalence obtained from the exact model when $SIS$ or $SIR$ propagation is considered on a complete graph, see \cite{armbruster2016elementaryIMA, armbruster2016elementaryJMB, armruster2016bounds} These result are valid for graphs with very simple structure, motivating research for finding upper and lower bounds for more complex models.

%A. Bátkai, I. Z. Kiss, E. Sikolya, P. L. Simon, Differential equation approximations of stochastic network processes: an operator semigroup approach, Netw. Heter. Media 7(2012), 43–58.
%D. Kunszenti-Kovács, P. L. Simon, Mean-field approximation of counting processes from a differential equation perspective, Electronic Journal of Qualitative Theory of Differential Equations 2016, No. 75, 1–17;
%B. Armbruster, Á. Besenyei, P. L. Simon, Bounds for the expected value of one-step processes, Commun. Math. Sci. 14(2016), No. 7, 1911–1923.
%B. Armbruster, E. Beck, An elementary proof of convergence to the mean-field equations for an epidemic model, IMA J. Appl. Math., published online on April 13, 2016.
%E. Beck and B. Armbruster. Elementary proof of convergence to the mean-field model for the SIR process. arXiv:1511.08572, 2015.

In this paper, we focus on the NIMFA model for $SIS$ epidemics and we will show that this model provides and upper bound on the exact model on arbitrary weighted and directed networks. In \cite{van2009virus} it is claimed that the NIMFA model overestimates the prevalence obtained from the exact system, however, the rigorous proof is not presented there. Here, this is done by using some well-known results from the theory of differential inequalities. The paper is structured as follows. In Section~\ref{sec:mod_form} the exact model is formulated and a bottom-up approach is used. This means that the model starts at the level of nodes and focuses on the probability of nodes being either susceptible or infected. Also here, the closure at the level of pairs is generalised beyond simply assuming that the state of nodes at the end of an edge are independent. In Section~\ref{sec:KamkeMuller} the main result is presented and it is here where we rigorously prove that NIMFA provides an upper bound on the exact model. In Section~\ref{sec:nonnegcorrel} we provide a proof, based on dynamical systems arguments, of the fact that epidemics are negatively correlated and point out that this is crucial for the proof of the main result. This is followed by the analysis on the closed model in Section~\ref{sec:anal_closed_model}. The paper concludes with a short discussion of the main findings and possible extensions.

%%%%%%%%%%%%%%%%%%%%%%%%%%%%%%%%%%%%%%%%%%%%%%%%%%%%%%%%%%%%%%%%%%%%%%%%%%%%%%%%%%%%%%%%%%
\section{Model formulation}\label{sec:mod_form}
%%%%%%%%%%%%%%%%%%%%%%%%%%%%%%%%%%%%%%%%%%%%%%%%%%%%%%%%%%%%%%%%%%%%%%%%%%%%%%%%%%%%%%%%%%

Consider a network with $N$ nodes and assume that no node has an edge to itself, but we allow for node $j$ to have an edge to node $i$ having some weight $g_{ij}$. Typically $g_{ij}=1$ if there is an edge from $j$ to $i$ and $0$ otherwise, but the model formulation works for any directed and weighted network and thus we can consider $g_{ij} \in [0,\infty)$ for $i,j=1, 2, \dots, N$.  We can use the adjacency matrix $G=(g_{ij})_{i,j=1, 2, \dots, N}$ to represent the network. We assume that the transmission rate from $j$ to $i$ is $\tau g_{ij}$. The recovery rate at each node is $\gamma$. The probability that node $i$ is infected at time $t$ is denoted by $\la I_i\ra(t)$. The aim is to derive exact and approximate differential equations for this function.

\subsection{Exact model}

It can be shown from first principles or from the exact master equations formulated in terms of the probabilities of all $2^N$ configurations~\cite{Sharkey2011} that $\la I_i\ra(t)$ satisfies the differential equation
\begin{equation}
\dot{\la I_i\ra}= \tau \sum_{j=1}^N g_{ij}\la S_iI_j\ra -\gamma\la I_i\ra,  \label{exacteq}
\end{equation}
where $\la S_iI_j\ra(t)$ is the probability that the pair consisting of node $i$ and node $j$ is of type $S-I$ at time $t$. This system is exact but not closed hence further differential equations or a closure is needed to determine the probability $\la I_i\ra(t)$. The differential equations for the pairs take the following form.
\begin{align}
\dot{\la S_iI_j\ra}&= \tau \sum_{k=1}^N g_{jk}\la S_iS_jI_k\ra- \tau \sum_{k=1}^N g_{ik}\la I_kS_iI_j\ra  -\tau g_{ij}\la S_iI_j\ra-\gamma\la S_iI_j\ra+\gamma \la I_iI_j \ra, \label{exactSI}\\
\dot{\la I_iS_j\ra}&= \tau \sum_{k=1}^N g_{ik}\la I_kS_iS_j\ra- \tau \sum_{k=1}^N g_{jk}\la I_iS_jI_k\ra  -\tau g_{ji}\la I_iS_j\ra-\gamma\la I_iS_j\ra+\gamma \la I_iI_j \ra, \label{exactSI}\\
\dot{\la I_iI_j\ra}&= \tau \sum_{k=1}^N g_{jk}\la I_iS_jI_k\ra +\tau \sum_{k=1}^N g_{ik}\la I_kS_iI_j\ra - 2\gamma \la I_iI_j\ra  +\tau g_{ij}\la S_iI_j\ra +\tau g_{ji}\la I_iS_j\ra, \label{exactII}\\
\dot{\la S_iS_j\ra}&=-\tau \sum_{k=1}^Ng_{ik}\la I_kS_iS_j\ra-\tau \sum_{k=1}^Ng_{jk}\la S_iS_jI_k\ra + \gamma\la S_iI_j\ra+\gamma\la I_iS_j\ra, \label{exactSS}
\end{align}
where $\la A_iB_jC_k\ra$ is the probability that the triple consisting of the nodes $i,j$ and $k$ is in the state $A-B-C$, and $(i,j)$ runs over all pairs satisfying $1\leq i<j\leq N$ and in all summations $k$ is different from $i$ and from $j$. Note that for any pair $(i,j)$ we have that the right hand sides of the four differential equations sums to 0, that is the sum $\la S_iI_j\ra+ \la I_iS_j\ra + \la I_iI_j\ra+ \la S_iS_j\ra$ remains constant in time. If this constant is 1 at the initial time instant, then it remains 1 for all time. %{\color{red} I am bit worried that this now does not apply to all directed and weighted networks or does it? Are we now assuming that $g_{ij}=g_{ji}$?}

\subsection{Closure at the level of pairs}

Closure at the level of pairs means that the joint probabilities $\la S_iI_j\ra$, $\la I_iS_j\ra$, $\la I_iI_j\ra$ and $\la S_iS_j\ra$ are expressed in terms of the marginal probabilities $\la I_i\ra$, $\la I_j\ra$, $\la I_i\ra$ and $\la S_j\ra$. This is always an approximation, as we will show below. However, we wish to rigorously quantify the accuracy of the approximation, or to be able to state whether a model based on closures gives upper or lower bounds for the exact values of the joint probabilities, or indeed if the prevalence from such a closed model is below or above the exact prevalence.

First we rigorously define what a closure relation means. The relation of the joint and marginal probabilities are shown in Table \ref{table1}, where
$a=\la I_iI_j\ra$, $b=\la S_iI_j\ra$, $c=\la I_iS_j\ra$ and $d=\la S_iS_j\ra$. All letters denote probabilities hence $a,b,c,d,p,q\in [0,1]$.
\begin{table}
\centering
            \begin{tabular}{| c  | c  |c  | c |}
                \hline
                  & $\la I_i\ra$ & $\la S_i\ra$ & \\
                 \hline
                $\la I_j\ra$  & a & b & q \\
                \hline
                $\la S_j\ra$  & c & d & 1-q \\
                \hline
                 & p & 1-p &  \\
                \hline
            \end{tabular}
\caption{The relation of the joint and marginal probabilities.} \label{table1}
\end{table}
The marginals can be expressed in terms of the joint probabilities as
\begin{equation}
a+b=q, \quad a+c=p, \quad c+d=1-q, \quad b+d=1-p. \label{joint}
\end{equation}
It can be immediately seen that $a+b+c+d=1$, hence the four equations are not independent, therefore they cannot be solved for the unknowns $a$, $b$, $c$ and $d$. This shows that the marginals do not determine the joint probabilities. However, once one of them is given then the remaining three are determined by system \eqref{joint}. A closure relation means that one of the joint probabilities is given by a certain algebraic relation involving the marginals, and the others are determined by system \eqref{joint}. Here we define the closure for the $II$ pairs, i.e. $a$ is specified in terms of $p$ and $q$. (One can equivalently express the probabilities of $SI$, $IS$ or $SS$ pairs.) So a closure will be a function $W:[0,1]\times [0,1] \to [0,1]$, for which
$$
\la I_iI_j\ra \approx  W(\la I_i\ra,\la I_j\ra) .
$$
In order to have a solution satisfying $a,b,c,d\in [0,1]$, the function $W$ must satisfy some conditions. On one hand, the first two equations of \eqref{joint} show that $a\leq p$ and $a\leq q$ must hold, that is $a\leq \min(p,q)$. On the other hand, the third and fourth equations should give a nonnegative value for $d$, hence $1-q-p+a\geq 0$ and $1-p-q+a\geq 0$ must hold, leading to $\max(p+q-1,0)\leq a$. (The $\max$ operation is needed because $p+q-1$ may be negative, when the lower bound for $a$ is zero.) Finally, it is natural to assume that $W$ is a symmetric function, i.e. $W(x,y)=W(y,x)$, leading to the following definition of the closure.
\begin{defi} \label{defi:closure}
A symmetric function $W:[0,1]\times [0,1] \to [0,1]$ is called a closure relation if
$$
\max(x+y-1,0)\leq W(x,y) \leq \min(x,y) \mbox{ for all } x,y\in [0,1].
$$
\end{defi}
We note that the most important closure relation is $W(x,y)=xy$ leading to the solution of system \eqref{joint} in the form
$$
a=pq, \quad b=(1-p)q, \quad c=p(1-q), \quad d=(1-p)(1-q) .
$$
This corresponds to the case when the states of node $i$ and $j$ are independent, i.e. the joint probability is the product of the marginals. To check that $W(x,y)=xy$ satisfies the inequalities in the definition is an easy exercise and is left to the Reader. One can also immediately see that $W(x,y)=\min(x,y)$ and $W(x,y)=\max(x+y-1,0)$ are proper closure relations.

Once a closure relation is chosen, the closed form of \eqref{exacteq} can be given as follows. Since $W(\la I_i\ra,\la I_j\ra)$ gives an approximation of $\la I_iI_j\ra$ and $\la S_iI_j\ra+ \la I_iI_j\ra= \la I_j\ra$, one can approximate $\la S_iI_j\ra$ as $\la I_j\ra - W(\la I_i\ra,\la I_j\ra)$. Hence the closed system is
\begin{equation}
\dot{X_i}= \tau \sum_{j=1}^N g_{ij}(X_j-W(X_i,X_j)) -\gamma X_i. \label{deqclosed}
\end{equation}
Solving this system for $X_i$ yields an approximation for $\la I_i\ra$.

The `art' of closing \eqref{exacteq} lies in choosing $W$ in a way in which $X_i$ is as close as possible to $\la I_i\ra$. In fact, the only function used in the literature is $W(x,y)=xy$ and the accuracy of the closure has been investigated only numerically. Here our goal is slightly different. We want to find upper and lower bounds to $\la I_i\ra$, i.e. to introduce closures $w$ and $W$ in such a way that for the corresponding solutions $x_i$ and $X_i$ of \eqref{deqclosed} the inequalities
$$
 x_i \leq \la I_i\ra \leq  X_i
$$
hold.

We will show that $W(x,y)=xy$ gives an upper bound and $W(x,y)=\min(x,y)$ gives a lower bound. However, improving these bounds by more sophisticated closures remains an open question. First, in the next section, we derive conditions on the closure $W$ ensuring that \eqref{deqclosed} gives an upper or lower bound.

\section{Bounds for the exact system} \label{sec:KamkeMuller}

The derivation of the upper bound is based on the fact that $S-I$ pairs are non-negatively, and $I-I$ and $S-S$ pairs are non-positively correlated \cite{cator2014nodal}, that is we have for any $i$ and $j$ that
\begin{equation}
\la S_iI_j\ra \leq \la S_i\ra \la I_j\ra, \quad \la I_iI_j\ra \geq \la I_i\ra \la I_j\ra , \quad \la S_iS_j\ra \geq \la S_i\ra \la S_j\ra.\label{SIcorrel}
\end{equation}
We will prove these inequalities in Section \ref{sec:nonnegcorrel} and use them to give an upper bound for the solution of system \eqref{exacteq}. The latter is done below. The first  relation leads to the differential inequality
\begin{equation}
\dot{\la I_i\ra}\leq \tau \sum_{j=1}^N g_{ij}\la S_i\ra \la I_j\ra -\gamma\la I_i\ra .  \label{diffineq}
\end{equation}
Based on this inequality let us introduce following system of differential equations, called individual-based model or N-intertwined mean-field approximation (NIMFA)~\cite{van2009virus,PietComputingSpringer}.
\begin{equation}
\dot{ Y_i}= \tau \sum_{j=1}^N g_{ij}(1- Y_i) Y_j -\gamma Y_i.  \label{nimfa}
\end{equation}
Since we have the same right hand side in \eqref{diffineq} and \eqref{nimfa} with inequality in the first, we expect that the NIMFA approximation yields an upper bound on the exact solution. We will prove this by using the comparison theory of ODEs based on the Kamka-M\"uller condition which is detailed below.

Consider the ODE $\dot x(t)=f(x(t))$ and the differential inequality $\dot y(t)\leq f(y(t))$ with a given differentiable function $f: \mathbb{R}^n \to \mathbb{R}^n$ subject to initial conditions satisfying $y(0)\leq x(0)$. The aim is to find conditions on $f$ ensuring that $y(t)\leq x(t)$ for $t\geq 0$. In what follows, the ordering relation for vectors is used in the following sense:
$$
u\leq v, \mbox{ if } u_i\leq v_i \,\,\, \forall~i=1,2,\ldots n, \qquad u< v, \mbox{ if } u_i< v_i \,\,\, \forall~i=1,2,\ldots n.
$$
A sufficient condition on $f$ for the desired inequality to hold is the Kamke-M\"uller condition \cite{kamke1932theorie,muller1927fundamentaltheorem}, which is equivalent to requiring that
$$
\mbox{ if } u\leq v \mbox{ and } u_i=v_i \Rightarrow f_i(u)\leq f_i(v)\,\,\, \forall~i=1,2,\ldots , n.
$$
This condition essentially means that the function in the $i$-th coordinate, $f_i$ is increasing in all coordinates $x_j$ for $j\neq i$.
However, this leads to an alternative sufficient condition, which if satisfied, allows us to call the $f$ function cooperative. More precisely, $f$ is called cooperative if
$$
\partial_j f_i \geq 0 \,\,\, \forall~i,j=1,2,\ldots n \,\,\, \mbox{and} \,\,\, i\neq j.
$$
It can be shown that if $f$ is cooperative in a convex domain, then it satisfies the Kamke-M\"uller condition. A detailed and comprehensive study of differential inequalities and comparison theory can be found in the book by Szarski \cite{szarski1965differential}. Cooperative systems generate monotone dynamical systems that are dealt with in the book chapter by Hirsch and Smith \cite{hirsch2005monotone} and the book by Smith \cite{smith2008monotone}. The main comparison result used here is the following.
\begin{lem}
Assume that $f$ satisfies the Kamke-M\"uller condition. If $\dot x(t)=f(x(t))$ and the differential inequality $\dot y(t)\leq f(y(t))$ holds and $y(0)\leq x(0)$,
then $y(t)\leq x(t)$ for $t\geq 0$. \label{lem:differentialinequality}
\end{lem}
Using this comparison result we can derive the following theorems about the upper and lower bounds.
\begin{theo} \label{theo:upperbound}
Consider a weighted and directed network $G$, the exact individual-based SIS model given by system \eqref{exacteq} and the closed system \eqref{deqclosed} with a closure $W$ satisfying (besides the conditions in Definition \ref{defi:closure})
\begin{enumerate}
\item[(i)] $y\mapsto y-W(x,y)$ is an increasing function,
\item[(ii)] $\la I_iI_j\ra \geq   W(\la I_i\ra,\la I_j\ra) $ for all $i,j$.
\end{enumerate}
Assuming that both models start with identical initial conditions, $\la I_i \ra(0)= X_i(0)$ ($i=1,2,\dots, N$), it follows that
 $\la I_{i} \ra(t) \leq X_{i}(t)$ $\forall~i=1,2, \dots, N$ and $\forall~t \geq 0$.
\end{theo}
%\begin{proof}
{\sc Proof}.
We start from the exact system,
\begin{align}
\dot{\la I_i\ra}&=\tau \sum_{j=1}^N g_{ij}\la S_iI_j\ra -\gamma\la I_i\ra = \tau \sum_{j=1}^N g_{ij}(\la I_j\ra-\la I_iI_j\ra) -\gamma\la I_i\ra \nonumber \\
&=\tau \sum_{j=1}^N g_{ij}(\la I_j\ra-W(\la I_i\ra,\la I_j\ra)) -\gamma\la I_i\ra+\tau \sum_{N=1}^{N}g_{ij}\left(W(\la I_i\ra,\la I_j\ra) - \la I_iI_j\ra\right)\nonumber \\
&\le\tau \sum_{j=1}^N g_{ij}(\la I_j\ra-W(\la I_i\ra,\la I_j\ra))-\gamma\la I_i\ra, \label{NonNegCorStrict}
\end{align}
where we have used (ii). The right-hand side of the closed system \eqref{deqclosed} can be given by the function $f: \mathbb{R}^N \to \mathbb{R}^N$ with coordinate functions
$$
f_i(x) = \tau \sum_{j=1}^N g_{ij} (x_j-W(x_i,x_j))- \gamma x_i .
$$
It can immediately be seen that the solution of the closed system satisfies the differential equation $\dot x(t)=f(x(t))$, and the solution of the exact system satisfies the differential inequality $\dot y(t)\leq f(y(t))$, with both systems starting from the same initial condition. According to (i) the coordinate function $f_i$ is increasing in the variable $x_j$, hence $f$ satisfies the Kamke-M\"uller condition. Therefore, the general comparison Lemma \ref{lem:differentialinequality} implies that $\langle I_i\rangle (t) \leq X_i(t)$ $\forall~i=1,2,\ldots, N$ and $\forall~t\geq 0$.

$\Box$
%\end{proof}

Now, using the closure $W(x,y)=xy$ we can prove that \eqref{nimfa} gives an upper bound. Namely, we have seen that $W(x,y)=xy$ satisfies the conditions in Definition \ref{defi:closure}. Moreover, $y-xy=y(1-x)$ is increasing in $y$ when $x\in [0,1]$, that is (i) holds. In Section \ref{sec:nonnegcorrel} we will prove that (ii) also holds, hence $W(x,y)=xy$ satisfies the conditions of Theorem \ref{theo:upperbound} leading to the following corollary.
\begin{cor}
Consider a weighted and directed network $G$, the exact individual-based SIS model given by system \eqref{exacteq} and the individual-based closed system \eqref{nimfa}.
Assuming that both models start with identical initial conditions, $\la I_i \ra(0)= Y_i (0)$ ($i=1,2,\dots, N$), it follows that
 $\la I_{i} \ra(t) \leq Y_{i}(t)$ $\forall~i=1,2, \dots, N$ and $\forall~t \geq 0$.
\end{cor}

Similarly to Theorem \ref{theo:upperbound}, the following result can be proved about the lower bound.

\begin{theo} \label{theo:lowerbound}
Consider a weighted and directed network $G$, the exact individual-based SIS model given by system \eqref{exacteq} and the closed system \eqref{deqclosed} with a closure $W$ satisfying (besides the conditions in Definition \ref{defi:closure})
\begin{enumerate}
\item[(i)] $y\mapsto y-W(x,y)$ is an increasing function,
\item[(ii)] $\la I_iI_j\ra \leq   W(\la I_i\ra,\la I_j\ra) $ for all $i,j$.
\end{enumerate}
Assuming that both models start with identical initial conditions, $\la I_i \ra(0)= X_i (0)$ ($i=1,2,\dots, N$), it follows that
 $\la I_{i} \ra(t) \geq X_{i}(t)$ $\forall~i=1,2, \dots, N$ and $\forall~t \geq 0$.
\end{theo}
It is easy to see that $W(x,y)=\min(x,y)$ satisfies the assumptions of the theorem. Namely, we have seen that $W(x,y)=\min(x,y)$ satisfies the conditions in Definition \ref{defi:closure}. Moreover, $y-\min(x,y)=\max(y-x,0)=\frac12 (|y-x|+y-x)$ is increasing in $y$, that is (i) holds. In the previous section we saw that $\la I_iI_j\ra \leq \la I_i\ra$ and $\la I_iI_j\ra \leq \la I_j\ra$, hence (ii) also holds. Thus $W(x,y)=\min(x,y)$ satisfies the assumptions of Theorem \ref{theo:lowerbound} leading to the following corollary.
\begin{cor}
Consider a weighted and directed network $G$, the exact individual-based SIS model given by system \eqref{exacteq} and the individual-based closed system \begin{equation}
\dot{X_i}= \tau \sum_{j=1}^N g_{ij}(X_j-\min(X_i,X_j)) -\gamma X_i. \label{deqmin}
\end{equation}
Assuming that both models start with identical initial conditions, $\la I_i \ra(0)= X_i(0)$ ($i=1,2,\dots, N$), it follows that
 $\la I_{i} \ra(t) \geq X_{i}(t)$ $\forall~i=1,2, \dots, N$ and $\forall~t \geq 0$.
\end{cor}

\section{$S-I$ pairs are non-negatively correlated} \label{sec:nonnegcorrel}

In this section we prove that $S-I$ pairs remain non-negatively correlated in the exact system \eqref{exacteq}-\eqref{exactSS} if they they are non-negatively correlated initially.
\begin{theo}
Let $\la S_i\ra$, $\la I_j\ra$ and $\la S_iI_j\ra$ solve system \eqref{exacteq}-\eqref{exactSS}. If $\la S_i\ra(0) \la I_j\ra(0)  -\la S_iI_j\ra(0) \geq 0$ for all $1\leq i<j\leq N$, then $\la S_i\ra(t) \la I_j\ra(t)  -\la S_iI_j\ra(t) \geq 0$ holds for all $t>0$ as well. \label{theo:SIcorrel}
\end{theo}
The proof of the theorem will be divided into several propositions.
\begin{prop}
$$\la S_i\ra \la I_j\ra  -\la S_iI_j\ra =  \la I_i I_j\ra \la S_i S_j\ra  -\la S_iI_j\ra \la I_iS_j\ra$$
\label{prop:correl}
\end{prop}
{\sc Proof}.
We will use the identities
$$
\la I_i I_j\ra+ \la S_i S_j\ra +\la S_iI_j\ra + \la I_iS_j\ra =1
$$
and
$$
\la S_i\ra = \la S_iI_j\ra + \la S_i S_j\ra , \quad  \la I_j\ra = \la S_iI_j\ra + \la I_i I_j\ra.
$$
Using these we obtain
$$
\la S_iI_j\ra\la I_i I_j\ra+ \la S_iI_j\ra\la S_i S_j\ra +\la S_iI_j\ra\la S_iI_j\ra + \la S_iI_j\ra\la I_iS_j\ra =\la S_iI_j\ra
$$
and
$$
\la S_i\ra \la I_j\ra = \la S_iI_j\ra\la S_i I_j\ra + \la S_iI_j\ra\la I_i I_j\ra + \la S_iS_j\ra\la S_iI_j\ra +\la S_iS_j\ra \la I_iI_j\ra .
$$
Taking the difference of the last two equations yields the desired relation.

$\Box$

Theorem \ref{theo:SIcorrel} will be proved by induction according to the the number of nodes in the network. First, we prove the theorem for a single edge, i.e. for $N=2$, for which system \eqref{exacteq}-\eqref{exactSS} takes the form
\begin{align}
\dot{\la I_1\ra}&= \tau \la SI\ra -\gamma\la I_1\ra, \label{exacteq1N2}\\
\dot{\la I_2\ra}&= \tau \la IS\ra -\gamma\la I_2\ra, \label{exacteq2N2}\\
\dot{\la SI\ra} &= -\tau \la SI\ra-\gamma\la SI\ra+\gamma \la II\ra, \label{exactSIN2}\\
\dot{\la IS\ra} &= -\tau \la IS\ra-\gamma\la IS\ra+\gamma \la II\ra, \label{exactISN2}\\
\dot{\la II\ra}&= - 2\gamma \la II\ra  +\tau \la SI\ra +\tau \la IS\ra, \label{exactIIN2}\\
\dot{\la SS\ra}&= \gamma \la SI\ra +\gamma \la IS\ra, \label{exactSSN2}
\end{align}
where we wrote $\la XY \ra$ instead of $\la X_1Y_2\ra$ for ease of notation.
\begin{prop}
The statement of Theorem \ref{theo:SIcorrel} holds for a graph with two connected nodes.
\label{prop:theoN2}
\end{prop}
{\sc Proof}.
We prove that
$$
A(t)=\la II\ra(t) \la SS\ra(t)  -\la SI\ra(t) \la IS\ra(t)
$$
is nonnegative if $A(0)\geq 0$, which proves the statement according to Proposition \ref{prop:correl}.

Using the differential equations \eqref{exacteq}-\eqref{exactSS} we can derive a differential equation for the function $A$ as follows.
\begin{align*}
\dot{A} &= \dot{\la II\ra} \la SS\ra + \la II\ra \dot{\la SS\ra} -\dot{\la SI\ra} \la IS\ra - \la SI\ra\dot{\la IS\ra} \\
&= \tau \la SI\ra \la SS\ra + \tau \la IS\ra \la SS\ra - 2\gamma \la II\ra \la SS\ra + \gamma \la SI\ra \la II\ra +\gamma \la IS\ra \la II\ra \\
& +\tau \la SI\ra\la IS\ra+\gamma\la SI\ra \la IS\ra-\gamma \la II\ra \la IS\ra +\tau \la IS\ra \la SI\ra+\gamma\la IS\ra \la SI\ra-\gamma \la II\ra \la SI\ra\\
&= \tau (\la SI\ra \la SS\ra  + \la IS\ra \la SS\ra +2\la SI\ra\la IS\ra) -2\gamma (\la II\ra \la SS\ra - \la IS\ra \la SI\ra)\\
&= -2(\tau+\gamma)A + b,
\end{align*}
where
$$
b= \tau(\la SI\ra \la SS\ra  + \la IS\ra \la SS\ra +2\la II\ra\la SS\ra ) .
$$
Thus $A$ satisfies an inhomogeneous linear differential equation. Multiplying this differential equation with $\exp(-2(\tau+\gamma)t)$ and integrating from $0$ to $t$ yields
$$
A(t)\mbox{e}^{-2(\tau+\gamma)t} - A(0) = \int\limits_0^t b(s) \mbox{e}^{-2(\tau+\gamma)s} \mbox{d}s .
$$
The non-negativity of $b(s)$ and $A(0)$ yields that $A(t)\geq 0$ for all nonnegative $t$.

$\Box$

Before proving the theorem in the general case, let us introduce the following notations,
\begin{equation}
A_{ij}(t)=\la S_i\ra(t) \la I_j\ra(t)  -\la S_iI_j\ra(t)=\la I_iI_j\ra(t) \la S_iS_j\ra(t)  -\la S_iI_j\ra(t) \la I_iS_j\ra(t), \label{Aij}
\end{equation}
where we have used Proposition~\ref{prop:correl}, and
\begin{equation}
P(S_i | S_k) = \frac{\la S_iS_k\ra}{\la S_k\ra},\label{PSiSk}
\end{equation}
which is the conditional probability that node $i$ is susceptible given that node $k$ is susceptible. The correlation of the conditional probabilities are given as
\begin{equation}
A_{ij}^k(t)=P(S_i | S_k) P(I_j | S_k) -P(S_i I_j| S_k),\label{Aijk}
\end{equation}
where the last term is the conditional probability that node $i$ is susceptible and node $j$ is infected given that node $k$ is susceptible. It is important to note that if
$A_{ij}^k(t)$ and $A_{ij}(t)$ identical as long as node $k$ is susceptible.

The proof of the theorem is based on the differential equation of $A_{ij}$ that is derived in the proposition below.
\begin{prop}
The functions $A_{ij}$ for $1\leq i<j\leq N$ satisfy the following system of inhomogeneous linear differential equations.
\begin{equation}
\dot{A_{ij}} = -b_{ij}A_{ij} + \tau \sum_{k=1}^N (g_{jk}P(S_j | S_k) A_{ki} + g_{ik}P(S_i | S_k) A_{kj}) + R_{ij},
\label{deAij}
\end{equation}
where
$$
b_{ij} = 2\gamma +\tau (g_{ij}+g_{ji}) + \tau \sum_{k=1}^N (g_{jk}+ g_{ik})
$$
and
$$
R_{ij} = \tau \la S_iS_j\ra (g_{ij} \la I_j\ra +g_{ji} \la I_i\ra) + \tau \sum_{k=1}^N \la S_k\ra(g_{jk} A_{ji}^k + g_{ik} A_{ij}^k) .
$$
\label{prop:deAij}
\end{prop}
{\sc Proof}.
For the ease of notation, we will write $\la XY \ra$ instead of $\la X_iY_j\ra$ and $\la I_kXY \ra$ instead of $\la I_kX_iY_j\ra$ and similarly for $\la X_iY_jI_k\ra$, where the indices $i$ and $j$ are fixed throughout the proof.

Differentiating \eqref{Aij} and using the differential equations \eqref{exacteq}-\eqref{exactSS} yield
\begin{equation}
\dot{A_{ij}} = \dot{\la II\ra} \la SS\ra + \la II\ra \dot{\la SS\ra} -\dot{\la SI\ra} \la IS\ra - \la SI\ra\dot{\la IS\ra} = Q_1+Q_2+ \tau \sum_{k=1}^N (g_{jk} Q_{3k} +g_{ik}Q_{4k}),\label{dotAij}
\end{equation}
where
\begin{align*}
Q_1&= \gamma (- 2 \la II\ra \la SS\ra + \la SI\ra \la II\ra +\la IS\ra \la II\ra +\la SI\ra \la IS\ra- \la II\ra \la IS\ra+\la IS\ra \la SI\ra- \la II\ra \la SI\ra ),\\
Q_2&= \tau g_{ij}\la SI\ra \la SS\ra + \tau g_{ji}\la IS\ra \la SS\ra +\tau g_{ij}\la SI\ra\la IS\ra +\tau g_{ji}\la IS\ra \la SI\ra \\
Q_{3k}&= \la ISI_k\ra \la SS\ra - \la SSI_k\ra \la II\ra - \la SSI_k\ra \la IS\ra + \la ISI_k\ra \la SI\ra \\
Q_{4k}&= \la I_kSI\ra \la SS\ra - \la I_kSS\ra \la II\ra + \la I_kSI\ra \la IS\ra - \la I_kSS\ra \la SI\ra .
\end{align*}

Each term will be simplified separately. Simple algebra leads to
$$
Q_1= -2\gamma (\la II\ra \la SS\ra - \la IS\ra \la SI\ra) =  -2\gamma A_{ij} .
$$
The expression for $Q_2$ can be reduced as follows
\begin{align*}
Q_2&= \tau g_{ij}(\la SI\ra \la SS\ra  + \la II\ra \la SS\ra -A_{ij}) + \tau g_{ji}(\la IS\ra \la SS\ra  + \la II\ra \la SS\ra -A_{ij}) \\
&= -\tau A_{ij} (g_{ij}+g_{ji}) + \tau g_{ij} \la SS\ra \la I_j\ra + \tau g_{ji} \la SS\ra \la I_i\ra ,
\end{align*}
where in the last step the identities $\la II\ra +\la SI\ra=\la I_j\ra$ and $\la II\ra +\la IS\ra=\la I_i\ra$ were used.

Based on the identities
$$
\la ISI_k\ra + \la ISS_k\ra=\la IS\ra , \quad \la SSI_k\ra + \la SSS_k\ra=\la SS\ra ,
$$
the term $Q_{3k}$ can be rewritten as
\begin{align*}
Q_{3k}&= \la ISI_k\ra \la SS\ra - \la SSI_k\ra \la IS\ra  + \la ISI_k\ra \la SI\ra - \la SSI_k\ra \la II\ra = (\la IS\ra- \la ISS_k\ra) \la SS\ra \\
&- (\la SS\ra- \la SSS_k\ra) \la IS\ra + (\la IS\ra-\la ISS_k\ra) \la SI\ra - (\la SS\ra -\la SSS_k\ra) \la II\ra\\
&= \la SSS_k\ra (\la IS\ra +\la II\ra) - \la ISS_k\ra (\la SS\ra + \la SI\ra) -  A_{ij} \\
&= (\la S_jS_k\ra -\la ISS_k\ra) \la I_i\ra - \la ISS_k\ra \la S_i\ra -  A_{ij} \\
&= \la S_jS_k\ra \la I_i\ra - \la ISS_k\ra -  A_{ij}.
\end{align*}
The transformations of the term $Q_{3k}$ can be completed by using \eqref{PSiSk} and \eqref{Aijk} as follows.
\begin{align*}
Q_{3k}&= \la S_jS_k\ra \la I_i\ra - \la ISS_k\ra -  A_{ij}\\
&= \la S_k\ra P(S_j | S_k) \la I_i\ra - \la S_k\ra P(I_iS_j | S_k)  -  A_{ij}\\
&= \la S_k\ra P(S_j | S_k) P(I_i | S_k) + P(S_j | S_k) (\la S_k\ra \la I_i\ra - \la I_iS_k\ra)  - \la S_k\ra P(I_iS_j | S_k)  -  A_{ij}\\
&= \la S_k\ra (P(S_j | S_k) P(I_i | S_k) - P(I_iS_j | S_k)) + P(S_j | S_k) A_{ki} -  A_{ij}\\
&= \la S_k\ra A_{ji}^k + P(S_j | S_k) A_{ki} -  A_{ij} ,
\end{align*}
where we used Proposition \ref{prop:correl} for $A_{ki}$.

Similar transformations lead to
$$
Q_{4k} = \la S_k\ra A_{ij}^k + P(S_i | S_k) A_{kj} -  A_{ij} .
$$
Substituting the expressions obtained for $Q_1$, $Q_2$, $Q_{3k}$ and $Q_{4k}$ into \eqref{dotAij} we get the desired differential equation for $A_{ij}$.

$\Box$

Now we are ready to prove Theorem \ref{theo:SIcorrel}.

{\sc Proof of Theorem \ref{theo:SIcorrel}}.
We prove that $A_{ij}(t)$ is nonnegative if $A_{ij}(0)\geq 0$, which proves the statement according to Proposition \ref{prop:correl}. This will be proved by induction according to $N$. The statement  for $N=2$ is proved in Proposition \ref{prop:theoN2}. Assuming that the statement is true for $N-1$ means that $A_{ij}^k\geq 0$, where node $k$ can be regarded as a newly added susceptible node
that does not invalidate what we already know. This effectively allows us to extend the model from $N-1$ to $N$ nodes, and thus completing the induction step.
%  because the condition that a given node is susceptible yields that we are concerning a graph with only $N-1$ nodes. \textbf{Eleg ennyi vagy reszletesebb magyarazat kellene?}
Let us now consider system \eqref{deAij} for all $i$ and $j$. This is an inhomogeneous system of linear differential equations in the form
$$
\dot x(t) = M x(t) + g(t),
$$
where $x(t)$ represents the vector containing all $A_{ij}$, $M$ is a cooperative matrix, i.e. its off-diagonal entries are nonnegative, and $g(t)\geq 0$. We prove that $x(0)\geq 0$ (coordinate-wise) implies $x(t)\geq 0$ for all $t\geq 0$. Namely, we have $\dot x(t) \geq M x(t)$, since $g(t)$ is non-negative. Consider the homogeneous system $\dot y(t) = M y(t)$ with initial condition $y(0)=0$. The solution is obviously $y(t)=0$ for all $t$. On the other hand, this system satisfies the Kamke-M\"uller condition, hence Lemma \ref{lem:differentialinequality} yields that $x(t)\geq y(t)$ for all $t\geq 0$, leading to the desired statement.

$\Box$

We note that the statement of the Theorem can be formulated in terms of the $I-I$ pairs. Using that $\la S_i\ra = 1- \la I_i\ra $ and $\la S_i I_j\ra = \la I_j\ra -\la I_i I_j\ra$, leads to the following.
\begin{cor}
Under the assumptions of Theorem \ref{theo:SIcorrel} the inequality $\la I_i\ra(t) \la I_j\ra(t)  \leq \la I_iI_j\ra(t)$ holds for all $t>0$. \label{cor:IIcorrel}
\end{cor}

%%%%%%%%%%%%%%%%%%%%%%%%%%%%%%%%%%%%%%%%%%%%%%%%%%%%%%%%%%%%%%%%%%%%%%%%%%%h%%%%%%%%%%%%%%%
\section{Analysis of the closed model}\label{sec:anal_closed_model}
%%%%%%%%%%%%%%%%%%%%%%%%%%%%%%%%%%%%%%%%%%%%%%%%%%%%%%%%%%%%%%%%%%%%%%%%%%%%%%%%%%%%%%%%%%

The goal in this section is to analyse system \eqref{deqclosed} from the dynamical system point of view. The closure $W$ satisfies the conditions in Definition \ref{defi:closure}, moreover, according to Corollary \ref{cor:IIcorrel} it satisfies $W(x,y)\geq xy$ as well. Thus our aim here is to understand the dynamical behaviour of system \eqref{deqclosed} when $W$ satisfies
\begin{equation}
xy\leq W(x,y) \leq \min(x,y) \mbox{ for all } x,y\in [0,1].\label{W_ini_cond}
\end{equation}
The lower bound $W(x,y)=xy$ is covered in \cite{lajmanovich1976deterministic}, and this result is presented first.

\begin{theo}\label{theo:LajmanYorke}
Given a directed, weighted, and strongly-connected network, $G$, let $\Lambda_{\text{max}}(G)$ be the largest eigenvalue of the adjacency matrix of $G$. Let the closure in \eqref{deqclosed} be given as $W(x,y)=xy$.
If $\gamma<\tau \Lambda_{\text{max}}(G)$, then a unique endemic (nonzero) steady state exists, and it is stable. Moreover, all of its coordinates are positive. If $\gamma>\tau \Lambda_{\text{max}}(G)$, then there is no endemic steady state and the disease-free steady state is stable.
\end{theo}

Concerning the upper bound $W(x,y)=\min(x,y)$, the following result can be easily proved.

\begin{prop}
Let $G$ be a directed and weighted network and let the closure in~\eqref{deqclosed} be given as $W(x,y)=\min(x,y)$. Then the only steady state is the disease-free and it is stable.
\end{prop}

Thus a transcritical bifurcation occurs when the closure is $W(x,y)=xy$ and there is no bifurcation when $W(x,y)=\min(x,y)$; that is the threshold behaviour disappears when such a crude closure is used. This however raises the question of studying the intermediate regime when $W$ is between the two extremes. Below we give a sufficient condition on closures to ensure that the threshold behaviour is maintained. 

We will consider closures where $W$ satisfies
\begin{equation}
xy\le W(x,y)\leq xy +V(x,y)\min(x,y)  \mbox{ for all } x,y\in [0,1],  \label{Wcond}
\end{equation}
where $V:[0,1]^2 \to [0,r]$, with some $r\in (0,1)$, is a continuous function satisfying $V(0,0)=0$ and $xy +V(x,y)\min(x,y) \le \min(x,y)$.
We note that the inequalities in~\eqref{Wcond} yield a sufficient condition for the existence of the transcritical bifurcation. This means that
it may be possible to identify closures that lead to transcritical bifurcation but do not satisfy condition~\eqref{Wcond}. 

We note that $W(x,y)=xy$ obviously satisfies this condition with $V(x,y)=0$, and a non-trivial example is $W(x,y)=\sqrt{xy}\min(\sqrt{x},\sqrt{y})$. For the latter, simple calculation shows that
there exists a $V(x,y)$ such that this is positive and bounded by a constant $r<1$.  Below we will prove that for any choice of $W$ that satisfies condition~\eqref{Wcond} the same threshold as in Theorem \ref{theo:LajmanYorke} is obtained.

The non-trivial steady state $x\in (0,1]^N$ of system \eqref{deqclosed} is given by
$$
\gamma x_i = \tau \sum_{j=1}^N g_{ij}(x_j-W(x_i,x_j))
$$
that will be rearranged using $x_j-W(x_i,x_j)=x_j-x_ix_j+x_ix_j-W(x_i,x_j)$ as
\begin{equation}
x_i(\alpha + (Gx)_i) = (Gx)_i -F_i(x) ,  \label{sseq}
\end{equation}
where $\alpha = \gamma /\tau$, $(Gx)_i$ is the $i$-th coordinate of the vector $Gx$ and
$$
F_i(x)=\sum_{j=1}^N g_{ij}(W(x_i,x_j)-x_ix_j) .
$$
Expressing $x_i$ from \eqref{sseq} we get the fixed point equation $x=T(x)$ for the non-trivial steady state with
\begin{equation}
T_i(x)=\frac{(Gx)_i -F_i(x)}{\alpha + (Gx)_i} . \label{Tix}
\end{equation}
We can immediately see that $T$ maps the unit cube $[0,1]^N$ into itself and the origin is its fixed point, representing the disease-free steady state. We will show that in the case $\gamma<\tau \Lambda$, that is $\alpha < \Lambda$, $T$ has a nontrivial fixed point in the interior of the cube, representing an endemic steady state. (Here $\Lambda=\Lambda_{\text{max}}(G)$ is the largest eigenvalue of the adjacency matrix of $G$.) The existence of this fixed point will be verified by using Brouwer's fixed point theorem on a suitably chosen convex subdomain of the cube not containing the origin. In order to achieve this goal we will need a few auxiliary results.

\begin{prop}
For a directed, weighted, and strongly-connected network, given by its adjacency matrix $G$, there exists a positive number $\mu$, for which the following holds. If $(Gx)_i< \eta$ and $x_i\geq 0$ for all $i=1,2, \ldots N$, then $|x|<\mu \eta$. \label{prop:Geta}
\end{prop}
{\sc Proof}.
Since the network is strongly connected every column of $G$ contains at least one nonzero entry. Hence $Gx\neq 0$ once $x_i\geq 0$ for all $i=1,2, \ldots N$ and $x\neq 0$. Therefore
$$
m=\min\{ |Gx|: \ x_i\geq 0, \ |x|=1  \} > 0.
$$
If $(Gx)_i< \eta$ for all $i=1,2, \ldots N$, then $|Gx| < \eta \sqrt{N}$. On the other hand, $|Gx| = |x| \left|G\frac{x}{|x|}\right| > m |x|$, implying $m|x| < \eta \sqrt{N}$. Hence the statement holds with $\mu = \sqrt{N}/m$.

$\Box$

For the next proposition we introduce a function, $h:[0,+\infty ) \to [0,1)$, which is defined by
$$
h(z)= \frac{z}{\alpha +z} .
$$
\begin{prop}
For any $\beta < 1/\alpha$ there is a $\omega >0$, such that $h(z)<\omega$ implies $\beta z < h(z)$, when $z\neq 0$. \label{prop:h}
\end{prop}
{\sc Proof}.
One can easily check that $h(0)=0$, $h'(0)=1/\alpha$, $h'(z)>0$ and $h''(z)<0$ for all $z\geq 0$, i.e. $h$ is increasing and concave. Hence a line with slope $\beta <h'(0)$ passing through the origin, lies below the graph of $h$ in a sufficiently short interval $(0,z_0)$. Then the statement holds with $\omega = h(z_0)$.

$\Box$

\begin{prop}
For any $x\in [0,1]^N$ we have $T_i(x)\geq (1-r) h((Gx)_i)$. \label{prop:Tih}
\end{prop}
{\sc Proof}.
According to \eqref{Wcond} we have $W(x_i,x_j)\leq x_ix_j +r x_j$, hence
$$
\sum_{j=1}^N g_{ij}(r x_j +x_ix_j - W(x_i,x_j)) \geq 0
$$
yielding $r(Gx)_i \geq F_i(x)$, that is $(Gx)_i -F_i(x) \geq (1-r) (Gx)_i$. Therefore \eqref{Tix} leads to
$$
T_i(x)\geq \frac{(1-r)(Gx)_i}{\alpha + (Gx)_i} = (1-r) h((Gx)_i).
$$

$\Box$

\begin{prop}
Let $\alpha< \Lambda$. Then there is a $\rho >0$, such that $T_i(x)<\rho$ for all $i$ implies $(Gx)_i< \Lambda T_i(x)$, when $(Gx)_i\neq 0$. \label{prop:rhoTi}
\end{prop}
{\sc Proof}.
Choose a small positive $\varepsilon$, for which $\Lambda (1-\varepsilon) > \alpha$ and introduce $\beta = 1/\Lambda (1-\varepsilon)<1/\alpha$. Choose $\omega$ to $\beta$ according to Proposition \ref{prop:h}. Choose a positive $\delta$ to $\varepsilon$ according to the continuity of $V$ given in \eqref{Wcond}, that is $|V(x,y)|<\varepsilon$, when $|x|,|y|<\delta$. Finally, determine $\kappa>0$ from
$$
\mu \alpha \frac{\kappa}{1-\kappa} = \delta,
$$
where $\mu$ is given in Proposition \ref{prop:Geta}. Now we show that choosing a positive $\rho$ satisfying
$$
\rho < \omega (1-\varepsilon) , \quad \rho < \kappa (1-r)
$$
the statement holds.

Using Proposition \ref{prop:Tih}, we get for any $i$ that
$$
\kappa (1-r) > \rho > T_i(x) \geq (1-r) h((Gx)_i) .
$$
Simple algebra shows that $\kappa > h((Gx)_i)$ implies $ (Gx)_i < \alpha \kappa / (1-\kappa)$. Hence, according to Proposition \ref{prop:Geta}, we have $|x|< \mu \alpha \frac{\kappa}{1-\kappa} = \delta$. The choice of $\varepsilon$ ensures that $W(x_i,x_j)-x_ix_j\leq \varepsilon x_j$, hence
$$
F_i(x) \leq \sum_{j=1}^N g_{ij} \varepsilon x_j = \varepsilon (Gx)_i .
$$
Therefore \eqref{Tix} leads to
$$
T_i(x)\geq \frac{(1-\varepsilon)(Gx)_i}{\alpha + (Gx)_i} = (1-\varepsilon) h((Gx)_i).
$$
Now,
$$
\omega (1-\varepsilon)  > \rho > T_i(x) \geq (1-\varepsilon) h((Gx)_i) .
$$
According to Proposition \ref{prop:h}, we get that $h((Gx)_i) > \beta (Gx)_i$, hence
$$
T_i(x) \geq (1-\varepsilon) h((Gx)_i) > (1-\varepsilon) \beta (Gx)_i = \frac{1}{\Lambda} (Gx)_i,
$$
and this completes the proof.

$\Box$

Now we are ready to prove the existence of the endemic steady state and the presence of a transcritical bifurcation.

\begin{theo}\label{theo:steadystates}
Given a directed, weighted, and strongly-connected network, $G$, let $\Lambda$ be the largest eigenvalue of the adjacency matrix of $G$. Let the closure $W$ in \eqref{deqclosed} satisfy \eqref{Wcond}. If $\gamma>\tau \Lambda$ then the origin is the only steady state of the system. In the case $\gamma<\tau \Lambda$, an endemic (nonzero) steady state also exists.
\end{theo}

{\sc Proof}.
We first consider the $\gamma>\tau \Lambda$ case and take a steady state $x\in [0,1]^N$. According to \eqref{sseq} and using that $F_i$ is nonnegative we get
$$
\gamma x_i \leq \tau (1-x_i)(Gx)_i\leq \tau (Gx)_i .
$$
It is easy to see that for two vectors, $u$ and $v$ with nonnegative coordinates, the inequality $0\leq u_i\leq v_i$ for all $i$ implies $|u|\leq |v|$. Hence for any nonzero steady state $x\in [0,1]^N$ we have
$$
\gamma |x|\leq \tau |G(x)| \leq \tau \Lambda |x| < \gamma |x|,
$$
where we used that $\Lambda$ is the largest eigenvalue of $G$. Hence there is no endemic steady state. We note that this part of the proof only used the fact that $W(x,y)\geq xy$, condition \eqref{Wcond} has not been used.

Let us turn to the case $\gamma<\tau \Lambda$. As we mentioned above, we will prove the existence of the endemic steady state by applying Brouwer's fixed point theorem to the mapping $T$ given in \eqref{Tix} in a suitably chosen domain. The goal is to exclude the origin from this domain, hence we introduce a half-space $S=\{ x\in \mathbb{R}^N : \ \la x-a,u\ra \geq 0 \}$ with some appropriately chosen vectors $a,u\in \mathbb{R}^N$. Then our domain will be $\Omega = [0,1]^N \cap S$. In order to have a nonempty intersection we will choose $a$ from the cube $[0,1]^N$ and $u$ will be the unique positive eigenvector of $G$ corresponding to the maximal eigenvalue $\Lambda$, that is $Gu=\lambda u$.

In order to prove the invariance of the domain $\Omega$ it is useful to determine the intersection points of the hyperplane given by $\la x-a,u\ra =0$ and the coordinate axes. The intersection point on the $i$-th axis is at $ c_i=\la a,u\ra /u_i$. It is easy to see that a point $y\in [0,1]^N$ is in $\Omega$ if there is a coordinate $i$, for which $y_i\geq c_i$. Namely, if $y\in [0,1]^N$ is not in $\Omega$, then $\la y,u\ra < \la a,u\ra = c_iu_i$ for all $i$, hence $y_iu_i < c_i u_i$ leading to $y_i<c_i$ for all $i$. Now, choose $a\in [0,1]^N$ in such a way that for all $i$ we have $c_i <\rho$ given in
Proposition \ref{prop:rhoTi}, that is $\la a,u\ra /u_i < \rho$ for all $i$.

We will prove that $T$ maps $\Omega$ into itself. We know that $T$ maps to $[0,1]^N$, hence we only need to prove that $\la T(x)-a,u\ra \geq 0$ holds for any $x\in \Omega$. We have previously shown that, if there is an index $i$ for which $T_i(x) \geq \rho$, then $T(x)\in \Omega$. Hence we only need to consider the case when $T_i(x)< \rho$ for all $i$. In this case we can apply Proposition
\ref{prop:rhoTi} yielding $T_i(x)> (Gx)_i/\Lambda$. Hence $T_i(x)u_i> (Gx)_iu_i/\Lambda$, leading to
\begin{align*}
\la T(x),u\ra &= \sum_{i=1}^N T_i(x)u_i > \frac{1}{\Lambda} \sum_{i=1}^N (Gx)_iu_i = \frac{1}{\Lambda} \sum_{i=1}^N u_i \sum_{i=j}^N g_{ij} x_j  \\
&= \frac{1}{\Lambda} \sum_{i=j}^N  x_j \sum_{i=1}^N g_{ij} u_i = \frac{1}{\Lambda} \sum_{i=j}^N  x_j (Gu)_j = \frac{1}{\Lambda} \sum_{i=j}^N  x_j \Lambda u_j = \la x,u\ra \geq \la a,u\ra .
\end{align*}
Thus we proved $\la T(x)-a,u\ra \geq 0$, which proves that $T$ maps the convex, compact domain $\Omega$ into itself, hence by Brouwer's fixed point theorem it has a fixed point in $\Omega$,
which is a nontrivial steady state completing the proof of the theorem.

$\Box$

We note that the stability of the steady states was also determined in the special case $W(x,y)=xy$. Further conditions on the smoothness of $W$ would make it possible to generalise the stability
result of Theorem 3.8 in \cite{kiss2016mathematics} for different choices of $W$.

%%%%%%%%%%%%%%%%%%%%%%%%%%%%%%%%%%%%%%%%%%%%%%%%%%%%%%%%%%%%%%%%%%%%%%%%%%%h%%%%%%%%%%%%%%%
\section{Discussion}
%%%%%%%%%%%%%%%%%%%%%%%%%%%%%%%%%%%%%%%%%%%%%%%%%%%%%%%%%%%%%%%%%%%%%%%%%%%%%%%%%%%%%%%%%%
In Section \ref{sec:nonnegcorrel} we proved that the closure $W(x,y)=xy$ satisfies the second assumption of Theorem \ref{theo:upperbound}, that is this closure leads to an upper bound of the exact
model. Similarly, in Theorem \ref{theo:lowerbound} we have shown that $W(x,y)=\min(x,y)$ yields a lower bound of the exact model. However, further work could focus on finding more accurate
upper and lower bounds with the possibility of constructing a sequence of closures whose limit is closer in some sense to the exact model.

The analysis of the closed model was presented in Section \ref{sec:anal_closed_model}. However, the investigation of the local and global stability for a general closure relation is still an open
question. Moreover, we have shown that the qualitative behaviour of the closed system depends on the closure and can be significantly different: it may or may not lead to a transcritical bifurcation.
The question of whether a closed system shares the same qualitative features as the exact model is an important one, and ideally the behaviour of the two systems should remain the same.
Thus identifying the closure or closures which separate different regimes, those that conserve the qualitative behaviour of the exact system versus those that do not, remains an important
direction for further research. One possible step towards achieving this may be to find closures that delimit closed models that have different qualitative behaviour when compared to each other, without considering their relation to the exact model.

%%%%%%%%%%%%%%%%%%%%%%%%%%%%%%%%%%%%%%%%%%%%%%%%%%%%%%%%%%%%%%%%%%%%%%%%%%%%%%%%%%%%%%%%%%
\section*{Acknowledgements} P\'eter L. Simon acknowledges support from Hungarian Scientific Research Fund, OTKA, (grant no. 115926).

%%%%%%%%%%%%%%%%%%%%%%%%%%%%%%%%%%%%%%%%%%%%%%%%%%%%%%%%%%%%%%%%%%%%%%%%%%%%%%%%%%%%%%%%%%

\bibliography{book_references_third_iteration}{}
\bibliographystyle{abbrv}
%
%\bibliography{Dyn_netw_chap}{}
%\bibliographystyle{spmpsci}
%\begin{thebibliography}{99}

%\bibitem{cator2014nodal}
%Cator, E., Van Mieghem, P., Nodal infection in Markovian susceptible-infected-susceptible and susceptible-infected-removed epidemics on networks are non-negatively correlated, Physical Review E 89.5 (2014): 052802.
%
%\bibitem{hirsch2005monotone}
%Hirsch, M.W., Smith, H., Monotone dynamical systems. In: A. Canada, P. Drábek, A. Fonda (eds.) Handbook of Differential Equations: Ordinary Differential Equations, vol. 2, pp. 239--357. Elsevier BV Amsterdam (2005)

%\bibitem{kamke1932theorie}
%Kamke, E., Zur theorie der systeme gew\"ohnlicher differentialgleichungen. ii, Acta Mathematica 58(1), 57–85 (1932).

%\bibitem{Lajmanovich}
%Lajmanovich, A., Yorke, J.A.: A deterministic model for gonorrhea in a nonhomogeneous population. Mathematical Biosciences 28(3), 221–236 (1976)
%
%\bibitem{muller1927fundamentaltheorem}
%M\"uller, M., \"U ber das fundamentaltheorem in der theorie der gew\"ohnlichen differentialgleichungen, Mathematische Zeitschrift 26(1), 619–645 (1927).
%
%\bibitem{smith2008monotone}
%Smith, H.L., Monotone dynamical systems: an introduction to the theory of competitive and cooperative systems, vol. 41. American Mathematical Soc. (2008)
%
%\bibitem{szarski1965differential}
%Szarski, J., Differential inequalities. Instytut Matematyczny Polskiej Akademi Nauk (Warszawa) (1965)

%\bibitem{van2011n}
%Van Mieghem, P., The n-intertwined SIS epidemic network model. Computing 93(2-4), 147--169 (2011).

%\end{thebibliography}

\end{document}